\documentclass[a4paper, BCOR=0.0mm, DIV=calc]{scrartcl}
\usepackage[T1]{fontenc}
\usepackage[utf8]{inputenc}
\usepackage[ngerman]{babel}
\usepackage[osf,sc]{mathpazo}
\usepackage{textcomp}
\usepackage{microtype}
\usepackage{amsmath, amssymb, amsfonts}
\usepackage{tabularx}
\KOMAoptions{twoside=false, twocolumn=false, headinclude=false, footinclude=false, mpinclude=false, pagesize=auto}
\recalctypearea

\setkomafont{section}{\mdseries\scshape\Large}
\setkomafont{title}{\mdseries\scshape}

\begin{document}
\title{Ein neuer Beweis, dass die Newton'sche Entwicklung der Potenzen des Binoms auch für gebrochene Exponenten gilt\footnote{
Originaltitel: "`Nova demonstratio, quod evolutio potestatum binomii Newtoniana etiam pro exponentibus fractis valeat"', erstmals publiziert in "`\textit{Nova Acta Academiae Scientarum Imperialis Petropolitinae} 5, 1789, pp. 52-58"', Nachdruck in "`\textit{Opera Omnia}: Series 1, Volume 16, pp. 112 - 121"', Eneström-Nummer E637, übersetzt von: Alexander Aycock, Textsatz: Artur Diener, im Rahmen des Projektes "`Eulerkreis Mainz"'}}
\author{Leonhard Euler}
\date{}
\maketitle
\paragraph{§1}
Wann immer in den Anfängen der Analysis die Potenzen des Binoms entwickelt wurden, pflegte das tatsächlich durch Multiplikation gemacht zu werden, während das Binom sooft mit sich selbst multipliziert wurde, so oft natürlich, wie der Exponent Einheiten enthält, und daher hat Newton für die unbestimmte  Potenz $(1+x)^n$ diese Progression der Terme gefolgert
\[
	1 + \frac{n}{1}x + \frac{n(n-1)}{1\cdot 2}xx + \frac{n(n-1)(n-2)}{1\cdot 2\cdot 3}x^3 + \mathrm{etc}.
\]
Die Gültigkeit dieser ist also nur für die Fälle, in denen der Exponent $n$ eine ganze positive Zahl ist, für bewiesen zu halten. Dass aber der selbe Ausdruck auch mit der Wahrheit verträglich ist, wann immer der Exponent $n$ entweder eine gebrochene oder negative oder sogar transzendente Zahl ist, haben viele Mathematiker zu zeigen versucht, deren Beweise aber entweder zu abstrus sind oder auch zu weit hergeholt, als dass sie im Haus der Analysis Platz finden können. Ich für meine Person habe auch vor vielen Jahren eine solchen Beweis gegeben, welcher die ersten Anfänge kaum zu überragen schein mag; neulich bin ich aber, nach einen anderen hineingestolpert, der mir freilich scheint die Aufgabe völlig zu erledigen; ich bin mir also sicher, dass diesen hier erörtert zu haben, den Mathematikern nicht unangenehm sein wird.
\paragraph{§2}
Von welcher Form auch immer der Exponent $n$ war, so lässt sich sicher annehmen, dass die Potenz selbst immer in eine Form solcher Art entwickelt werden kann, dass
\[
	(1+x)^n = 1 + Ax + Bxx + Cx^3 + Dx^4 + \mathrm{etc}
\]
ist. Hier bezeichnen natürlich die Großbuchstaben gewisse durch die Exponenten $n$ zu bestimmende Zahlen, wo wir freilich schon wissen, sooft $n$ eine ganze Zahl positive Zahl war, dass dann
\[
	A = \frac{n}{1},\quad B = \frac{n-1}{2}A,\quad C = \frac{n-2}{3}B,\quad D = \frac{n-3}{4}C,\quad \mathrm{etc}
\]
sein wird; dennoch müssen diese Werte auch durch die Methode, die ich hier erörtern will, gefolgert werden, woher zugleich klar werden wird, dass dieselben sogar immer Geltung haben, wenn auch der Exponent $n$ keine ganze positive Zahl war.
\paragraph{§3}
Sofort wird es hier aber förderlich sein, bemerkt zu haben, dass der erste Term der angenommenen Reihe auf gewohnte Weise der Einheit gleich gesetzt wird, weil wir ja wissen, wenn $x=0$ war, in welchem Fall alle Terme, die dem ersten folgen, verschwinden, dass der Wert der Potenz $1^n$ gewiss immer gleich $1$ ist, welcher Wert auch immer für $n$ angenommen wurde. Darauf ist klar, dass im Fall, in dem $n=0$ ist, der Wert der Potenz $(1+x)^0$ immer der Einheit selbst gleich wird, weil sogar in den Prinzipen der Analysis hinreichend und darüber hinaus klar dargetan wurde, dass immer $z^0 = 1$ ist. Daher folgt also in dem Fall, in dem $n=0$ ist, dass auch die Werte aller Buchstaben $A$, $B$, $C$, $D$, etc verschwinden müssen, damit natürlich für den ganzen Ausdruck der Wert gleich $1$ hervorgeht. Deshalb ist es notwendig, dass die einzelnen dieser Buchstaben den Faktor $n$ involvieren, so wie es in den von Newton festgesetzen Werte passiert.
\paragraph{§4}
Wir wollen also nun den Exponenten unserer Potenz $n$ um die Einheit vermehren und auf ähnliche Weise
\[
	(1+x)^{n+1} = 1 + A'x + B'xx + C'x^3 + D'x^4 + \mathrm{etc}
\]
festsetzen, wo klar ist, dass diese mit dem Strich versehenen Buchstaben $A'$, $B'$, $C'$, etc aus den vorhergehenden entstehen müssen, wenn in diese überall $n+1$ anstelle von $n$ geschrieben wird. Weil aber
\[
	(1+x)^{n+1} = (1+x)(1+x)^n
\]
ist, ist klar, dass die neue Reihe aus der ersten entstehen muss, wenn diese mit $1+x$ multipliziert werden; dann wird sich aber das nach Potenzen von $x$ geordnete Produkt so verhalten
\begin{align*}
1 + A&x + Bxx + Cx^3 + Dx^4 + Ex^5 + \mathrm{etc} \\
    +&x + Axx + Bx^3 + Cx^4 + Dx^5 + \mathrm{etc},
\end{align*}
welche zwei Reihen also zusammengenommen unserer neuen Reihe gleich sein müssen.
\paragraph{§5}
Daher erhält man aber durch Vergleichen der einzelnen Terme miteinander folgende Gleichheiten:
\begin{alignat*}{2}
	A' = A + 1 &\quad \text{oder} \quad &A' - A = 1 \\
	B' = B + A &\quad \text{oder} \quad &B' - B = A \\
	C' = C + B &\quad \text{oder} \quad &C' - C = B \\
	D' = D + C &\quad \text{oder} \quad &D' - D = C \\
	E' = E + D &\quad \text{oder} \quad &E' - E = D \\
	etc
\end{alignat*} 
woher man einsieht, wie zwei gewisse folgenden Buchstaben auf den Vorhergehenden bezogen werden. Wenn z.\,B. im Allgemeinen der Buchstabe $M$ dem Buchstaben $N$ gleich wird, wird
\[
	N' - N = M
\]
sein; daher geht die ganze Aufgabe darauf zurück, wie, wenn der Wert des Buchstaben $N$ gefunden worden ist, der folgende Buchstabe untersucht werden muss, sodass, wenn $n+1$ bei ihm anstelle von $n$ geschrieben wird und der daher resultierende Wert mit $N'$ gekennzeichnet wird, $N' - N = M$ ist.
\paragraph{§6}
Diese Untersuchung also wollen wir auf folgende Weise unternehmen, indem wir von einfachsten Fällen aus anfangen, für welches Ziel wir die folgenden Lemmata vorausschicken wollen.
\section*{Lemma 1}
\paragraph{§7}
Wenn $N=\alpha n$ war, wird $N' = \alpha (n+1)$ sein und daher $N'-N = \alpha$; daher ist andererseits klar, wenn $M=\alpha$ war, dass dann $N = \alpha n$ sein wird. Hier könnte man freilich entgegnen, dass der Umkehrschluss nicht hinreichend gewiss ist. Wenn wir nämlich $N = \alpha n + \beta$ gesetzt hätten, wäre $N' = \alpha (n+1) + \beta$ hervorgegangen und daher $N' - N = \alpha$; daher scheint also gefolgert zu sein, wenn $M=\alpha$ war, dass dann auf allgemeine Weise $N = \alpha n + \beta$ gesetzt werden muss. Aber wir haben schon anfangs angemerkt, dass unsere Koeffizienten $A$, $B$, $C$, $D$, etc so beschaffen sein müssen, dass sie für $n=0$ gesetzt verschwinden; weil daher $N$ jeden einzigen dieser Buchstaben bezeichnet, ist es klar, dass notwendigerweise $\beta = 0$ genommen werden muss, und so ist vollkommen dargetan worden, sooft $M=\alpha$ war, dass $N = \alpha n$ genommen werden muss.
\section*{Lemma 2}
\paragraph{§8}
Wenn $N = \alpha n(n-1)$ war, wird $N' = \alpha (n+1)n$ sein, woher $N' - N = 2\alpha n$ entsteht, sodass in diesem Fall $M = 2\alpha n$ ist. Indem wir anstelle von $2\alpha$ also $\alpha$ schreiben, folgern wir, sooft $M = \alpha n$ war, dass dann gewiss $N = \frac{1}{2}\alpha n(n-1)$ sein wird; dieser Wert, weil er schon für $n=0$ gesetzt verschwindet, kann keinen konstanten Zuwachs erhalten, wie wir gerade zuvor gezeigt haben; das ist auch im Folgenden festzuhalten.
\section*{Lemma 3}
\paragraph{§9}
Wenn $N = \alpha n(n-1)(n-2)$ war, wird $N' = \alpha (n+1)n(n-1)$ sein und daher $N' - N = 3\alpha n(n-1)$, sodass in diesem Fall $M = 3\alpha n(n-1)$ ist. Daher schließen wir also andererseits, sooft $M = an(n-1)$ war, dass gewiss $N = \frac{1}{3}an(n-1)(n-2)$ sein wird.
\section*{Lemma 4}
\paragraph{§10}
Wenn $N = \alpha n(n-1)(n-2)(n-3)$ war, wird $N' = \alpha (n+1)n(n-1)(n-2)$ sein, der gemeinsame Faktor welcher beider Formeln $\alpha n(n-1)(n-2)$ ist, aus welchem $N' - N = \alpha n(n-1)(n-2)((n+1)-(n-3))$ sein wird und daher $N'-N = 4\alpha n(n-1)(n-2) = M$; und daher schlussfolgern wir andererseits, sooft $M = an(n-1)(n-2)$ war, dass gewiss $N = \frac{1}{4}an(n-1)(n-2)(n-3)$ sein wird.
\section*{Lemma 5}
\paragraph{§11}
Wenn $N = \alpha n(n-1)(n-2)(n-3)(n-4)$ war, wird $N' = \alpha (n+1)n(n-1)(n-2)(n-3)$ sein und daher berechnet man, dass $N' - N = M = 5\alpha n (n-1)(n-2)(n-3)$ sein wird, woraus wir andererseits schlussfolgern, sooft $M = an(n-1)(n-2)(n-3)$ war, dass dann immer $N = \frac{1}{5}an(n-1)(n-2)(n-3)(n-4)$ sein wird.
\paragraph{§12}
Daher durchschaut man schon hinreichend deutlich, wenn 
\[
	M = an(n-1)(n-2)(n-3)(n-4)
\]
ist, dass dann gewiss
\[
	N = \frac{1}{6}an(n-1)(n-2)(n-3)(n-4)(n-5)
\]
sein wird und daher im Allgemeinen, sooft
\[
	M = an(n-1)(n-2)(n-3)\cdots (n-\lambda),
\]
dass dann gewiss
\[
	N = \frac{1}{\lambda +2}an(n-1)(n-2)\cdots (n-\lambda -1)
\]
sein wird, welche allgemeine Form alle vorhergehenden Lemmata zugleich in sich einfasst.
\paragraph{§13}
Weil nun die Buchstaben $M$ und $N$ im Allgemeinen zwei in der Reihe der Buchstaben $A$, $B$, $C$, $D$, etc aufeinander folgende Terme bezeichnen, wollen wir hier die im Allgemeinen entwickelte Gleichung
\[
	N'-N = M
\]
auf die einzelnen oben in §$5$ gefundenen Gleichheiten anwenden; weil deren erste $A' - A = 1$ ist, wird hier $M=1$ sein, und daher berechnen wir aus Lemma $1$, dass $N$ oder $A=n$ sein wird; dieser Wert wird also gewiss wahr sein, welche Zahl auch immer für den Exponenten $n$ angenommen wird, weil ja die oberen Rechnungen nie auf ganze Zahlen beschränkt waren.
\paragraph{§14}
Wir wollen also zur zweiten Gleichung 
\[
	B' - B = A
\]
weiter gehen, wo, weil wir gerade gefunden haben, dass $n=A$ ist, aufgrund von Lemma $2$ $M=A=n$ sein wird, woher uns der Wert von $N$ hier 
\[
	B = \frac{1}{2}n(n-1)
\]
oder
\[
	B =\frac{n}{1}\frac{n-1}{2}
\]
liefern wird, völlig wie die Newton'sche Entwicklung es angegeben hatte.
\paragraph{§15}
Weil also weiter die dritte Gleichheit
\[
	C' - C = B
\]
war, wird Lemma $3$ zur Hilfe nehmend,
\[
	M = \frac{1}{2}n(n-1)
\]
sein und daher $a=\frac{1}{2}$, woher der dort für $N$ gefundene Wert uns hier
\[
	C = \frac{1}{6}n(n-1)(n-2)
\] 
geben wird oder
\[
	C = \frac{n}{1}\cdot \frac{n-1}{2}\cdot \frac{n-2}{3},
\]
genauso wie es die Newton'sche Entwicklung hält.
\paragraph{§16}
Unsere vierte Gleichheit war
\[
	D' - D = C,
\]
die mit Lemma $4$ verglichen
\[
	M = C = \frac{1}{6}n(n-1)(n-2)
\]
liefert, sodass $a = \frac{1}{6}$ ist; daher wird der Buchstabe $N$ uns also
\[
	D = \frac{1}{24}n(n-1)(n-2)(n-3)
\]
liefern oder
\[
	D = \frac{n}{1}\cdot \frac{n-1}{2}\cdot \frac{n-2}{3}\cdot \frac{n-3}{4}.
\]
\paragraph{§17}
Nach dem Fund dieses Wertes wollen wir zur fünften Gleichheit gelangen
\[
	E' - E = D,
\]
die uns für Lemma $5$
\[
	M = D = \frac{1}{24}n(n-1)(n-2)(n-3)
\]
beschafft, sodass hier $a=\frac{1}{24}$ ist; daher liefert der dort für $N$ angegebene Wert für den gegenwärtigen Fall
\[
	E = \frac{1}{120}n(n-1)(n-2)(n-3)(n-4)
\]
oder auf gewohnte Weise
\[
	E = \frac{n}{1}\cdot\frac{n-1}{2}\cdot\frac{n-2}{3}\cdot\frac{n-3}{4}\cdot\frac{n-4}{5}.
\]
\paragraph{§18}
Es wäre völlig überflüssig, diese Fälle weiter zu verfolgen, weil schon klarer als das Mittagslicht scheint, dass für die einzelnen folgenden Buchstaben natürlich notwendigerweise dieselben Werte hervorgehen müssen, welche die Newton'sche Entwicklung schon gelehrt hat, und dieser Beweis scheint für die Natur der Sache dermaßen geeignet, dass jenem auch in den ersten Elementen der Analysis ein Platz nicht versagt werden kann. Ja sogar die ganze Methode, die wir hier benutzt haben, behält die ganze Kraft bei, auch wenn sogar der Exponent $n$ als imaginär angesehen wird. 
\end{document}